\documentclass{gtart}

\def\ifplaintex{\expandafter\ifx\csname documentclass\endcsname\relax}

\def\gtp{{\mathsurround=0pt\it $\cal G\mskip-2mu$eometry \&\ 
$\cal T\!\!$opology $\cal P\!$ublications}}  

\def\recd{{\small Received:\qua\receiveddate\ifx\reviseddate\relax
\else\qquad Revised:\qua\reviseddate\fi\par}} 

\def\lognumber#1{\def\thelognumber{#1}}
\def\volumenumber#1{\def\thevolumenumber{#1}}
\def\volumeyear#1{\def\thevolumeyear{#1}}
\def\papernumber#1{\def\thepapernumber{#1}}
\def\pagenumbers#1#2{\def\startpage{#1}\def\finishpage{#2}}
\def\published#1{\def\publishdate{#1}}

\def\received#1{\def\receiveddate{#1}}
\def\revised#1{\def\reviseddate{#1}}
\def\accepted#1{\def\accepteddate{#1}}

\def\asciiaddress#1{\def\theasciiaddress{#1}}

\let\\\par\let\thelognumber\relax\let\thevolumenumber\relax
\let\thepapernumber\relax\let\thevolumeyear\relax\let\startpage\relax
\let\finishpage\relax\let\publishdate\relax\let\receiveddate\relax
\let\reviseddate\relax\let\accepteddate\relax\let\theasciititle\relax
\let\theasciiauthors\relax\let\theasciiaddress\relax
\let\theasciiabstract\relax

\let\theasciiemail\relax

\ifplaintex
\font\logobig=cmssbx10 scaled 3836
\font\logomed=cmssbx10 scaled 2557
\else
\font\logobig=cmssbx10 scaled 4200
\font\logomed=cmssbx10 scaled 2800
\fi

\long\def\makeagttitle{   
\count0=\startpage
\agt\hfill      
\hbox to 45truept{\vbox to 0pt{\vglue -13truept{\logomed A\kern -.37em{\logobig 
T}\kern -.38em G}\vss}\hss}
\break
{\small Volume \thevolumenumber\ (\thevolumeyear)
\startpage--\finishpage\nl
Published: \publishdate}

\vglue .25truein

{\parskip=0pt\leftskip 0pt plus
1fil\def\\{\par\smallskip}{\Large\bf\thetitle}\par\medskip} \vglue
0.05truein

{\parskip=0pt\leftskip 0pt plus 1fil\def\\{\par}{\sc\theauthors}
\par\medskip}%
 
\vglue 0.03truein 

{\small\leftskip 25truept\rightskip 25truept{\bf Abstract}\stdspace\theabstract

{\bf AMS Classification}\stdspace\theprimaryclass
\ifx\thesecondaryclass\relax\else; \thesecondaryclass\fi\par
{\bf Keywords}\stdspace \thekeywords\par}\vglue 7truept

}   

\ifplaintex
\font\phead=cmsl9 scaled 950
\font\pnum=cmbx10 scaled 913
\font\pfoot=cmsl9 scaled 950
\headline{\vbox to 0pt{\vskip -4.5mm\line{\small\phead\ifnum
\count0=\startpage ISSN 1472-2739 (on-line) 1472-2747 (printed)
\hfill {\pnum\folio}\else\ifodd\count0\def\\{ }%
\ifx\theshorttitle\relax\thetitle\else\theshorttitle\fi\hfill{\pnum\folio}
\else\def\\{ and }{\pnum\folio}\hfill\ifx\theshortauthors\relax\theauthors
\else\theshortauthors\fi\fi\fi}\vss}}
\footline{\vbox to 0pt{\vglue 0mm\line{\small\pfoot\ifnum\count0=\startpage
\copyright\ \gtp\hfill\else
\agt, Volume \thevolumenumber\ (\thevolumeyear)\hfill\fi}\vss}}
\else
\font=cmsl9 scaled 1050
\font\lnum=cmbx10 
\font=cmsl9 scaled 1050
\makeatletter
\def\@oddhead{{\small\lhead\ifnum\count0=\startpage ISSN 1472-2739 
(on-line) 1472-2747 (printed)\hfill {\lnum\number\count0}\else\ifodd\count0
\def\\{ }\ifx\theshorttitle\relax \thetitle \else\theshorttitle\fi\hfill
{\lnum\number\count0}\else\def\\{ and }{\lnum\number\count0}
\hfill\ifx\theshortauthors\relax 
\theauthors\else\theshortauthors\fi\fi\fi}}\def\@evenhead{\@oddhead}
\def\@oddfoot{\small\lfoot\ifnum\count0=\startpage\copyright\ \gtp\hfill\else
\agt, Volume \thevolumenumber\ (\thevolumeyear)\hfill\fi}
\def\@evenfoot{\@oddfoot}
\makeatother
\fi
\let\maketitlepage\makeagttitle
\let\makeshorttitle\maketitlepage
\let\maketitle\maketitlepage

\newwrite\gtoutfile
\long\gdef\makeheadfile{  
{\def\\{, }\def\s{ }
\immediate\openout\gtoutfile head.xxx
\immediate\write\gtoutfile{To: math@arxiv.org}
\immediate\write\gtoutfile{Subject: put OR rep NNNNN:ppppp}
\immediate\write\gtoutfile{--text follows this line--}
\immediate\write\gtoutfile{Proxy-for: \ifx\theasciiauthors\relax
\theauthors\else\theasciiauthors\fi\s<\ifx\theasciiemail\relax\theemail\else\theasciiemail\fi>}
\immediate\write\gtoutfile{\noexpand\\}
\immediate\write\gtoutfile{Authors: \ifx\theasciiauthors\relax
\theauthors\else\theasciiauthors\fi}
{\def\\{ }\immediate\write\gtoutfile{Title: \ifx\theasciititle\relax
\thetitle\else\theasciititle\fi}}
\immediate\write\gtoutfile{Subj-class: GT or SG, GR etc}
\immediate\write\gtoutfile{MSC-class: \theprimaryclass\ifx\thesecondaryclass\relax\else, \thesecondaryclass\fi}
\immediate\write\gtoutfile{Journal-ref: Algebr. Geom. Topol. \thevolumenumber\s
(\thevolumeyear) \startpage-\finishpage}
\immediate\write\gtoutfile{Comments: Published by Algebraic and
Geometric Topology at}
\immediate\write\gtoutfile{\s\s\s  http://www.maths.warwick.ac.uk/agt/AGTVol\thevolumenumber/agt-\thevolumenumber-\thepapernumber.abs.html}
\immediate\write\gtoutfile{\noexpand\\}
\immediate\write\gtoutfile{}
\ifx\theasciiabstract\relax
\immediate\write\gtoutfile{\theabstract}\else
\immediate\write\gtoutfile{\theasciiabstract}\fi
\immediate\write\gtoutfile{}
\immediate\write\gtoutfile{\noexpand\\}
\immediate\write\gtoutfile{}
\immediate\closeout\gtoutfile}}  

\def\maketitlepage{\makeagttitle\makeheadfile}
\let\makeshorttitle\maketitlepage
\let\maketitle\maketitlepage

\def\ifplaintex{\expandafter\ifx\csname documentclass\endcsname\relax}

\def\gtp{{\mathsurround=0pt\it $\cal G\mskip-2mu$eometry \&\ 
$\cal T\!\!$opology $\cal P\!$ublications}}  

\def\recd{{\small Received:\qua\receiveddate\ifx\reviseddate\relax
\else\qquad Revised:\qua\reviseddate\fi\par}} 

\def\lognumber#1{\def\thelognumber{#1}}
\def\volumenumber#1{\def\thevolumenumber{#1}}
\def\volumeyear#1{\def\thevolumeyear{#1}}
\def\papernumber#1{\def\thepapernumber{#1}}
\def\pagenumbers#1#2{\def\startpage{#1}\def\finishpage{#2}}
\def\published#1{\def\publishdate{#1}}

\def\received#1{\def\receiveddate{#1}}
\def\revised#1{\def\reviseddate{#1}}
\def\accepted#1{\def\accepteddate{#1}}

\def\asciiaddress#1{\def\theasciiaddress{#1}}

\let\\\par\let\thelognumber\relax\let\thevolumenumber\relax
\let\thepapernumber\relax\let\thevolumeyear\relax\let\startpage\relax
\let\finishpage\relax\let\publishdate\relax\let\receiveddate\relax
\let\reviseddate\relax\let\accepteddate\relax\let\theasciititle\relax
\let\theasciiauthors\relax\let\theasciiaddress\relax
\let\theasciiabstract\relax

\let\theasciiemail\relax

\ifplaintex
\font\logobig=cmssbx10 scaled 3836
\font\logomed=cmssbx10 scaled 2557
\else
\font\logobig=cmssbx10 scaled 4200
\font\logomed=cmssbx10 scaled 2800
\fi

\long\def\makeagttitle{   
\count0=\startpage
\agt\hfill      
\hbox to 45truept{\vbox to 0pt{\vglue -13truept{\logomed A\kern -.37em{\logobig 
T}\kern -.38em G}\vss}\hss}
\break
{\small Volume \thevolumenumber\ (\thevolumeyear)
\startpage--\finishpage\nl
Published: \publishdate}

\vglue .25truein

{\parskip=0pt\leftskip 0pt plus
1fil\def\\{\par\smallskip}{\Large\bf\thetitle}\par\medskip} \vglue
0.05truein

{\parskip=0pt\leftskip 0pt plus 1fil\def\\{\par}{\sc\theauthors}
\par\medskip}%
 
\vglue 0.03truein 

{\small\leftskip 25truept\rightskip 25truept{\bf Abstract}\stdspace\theabstract

{\bf AMS Classification}\stdspace\theprimaryclass
\ifx\thesecondaryclass\relax\else; \thesecondaryclass\fi\par
{\bf Keywords}\stdspace \thekeywords\par}\vglue 7truept

}   

\ifplaintex
\font\phead=cmsl9 scaled 950
\font\pnum=cmbx10 scaled 913
\font\pfoot=cmsl9 scaled 950
\headline{\vbox to 0pt{\vskip -4.5mm\line{\small\phead\ifnum
\count0=\startpage ISSN 1472-2739 (on-line) 1472-2747 (printed)
\hfill {\pnum\folio}\else\ifodd\count0\def\\{ }%
\ifx\theshorttitle\relax\thetitle\else\theshorttitle\fi\hfill{\pnum\folio}
\else\def\\{ and }{\pnum\folio}\hfill\ifx\theshortauthors\relax\theauthors
\else\theshortauthors\fi\fi\fi}\vss}}
\footline{\vbox to 0pt{\vglue 0mm\line{\small\pfoot\ifnum\count0=\startpage
\copyright\ \gtp\hfill\else
\agt, Volume \thevolumenumber\ (\thevolumeyear)\hfill\fi}\vss}}
\else
\font=cmsl9 scaled 1050
\font\lnum=cmbx10 
\font=cmsl9 scaled 1050
\makeatletter
\def\@oddhead{{\small\lhead\ifnum\count0=\startpage ISSN 1472-2739 
(on-line) 1472-2747 (printed)\hfill {\lnum\number\count0}\else\ifodd\count0
\def\\{ }\ifx\theshorttitle\relax \thetitle \else\theshorttitle\fi\hfill
{\lnum\number\count0}\else\def\\{ and }{\lnum\number\count0}
\hfill\ifx\theshortauthors\relax 
\theauthors\else\theshortauthors\fi\fi\fi}}\def\@evenhead{\@oddhead}
\def\@oddfoot{\small\lfoot\ifnum\count0=\startpage\copyright\ \gtp\hfill\else
\agt, Volume \thevolumenumber\ (\thevolumeyear)\hfill\fi}
\def\@evenfoot{\@oddfoot}
\makeatother
\fi
\let\maketitlepage\makeagttitle
\let\makeshorttitle\maketitlepage
\let\maketitle\maketitlepage

\newwrite\gtoutfile
\long\gdef\makeheadfile{  
{\def\\{, }\def\s{ }
\immediate\openout\gtoutfile head.xxx
\immediate\write\gtoutfile{To: math@arxiv.org}
\immediate\write\gtoutfile{Subject: put OR rep NNNNN:ppppp}
\immediate\write\gtoutfile{--text follows this line--}
\immediate\write\gtoutfile{Proxy-for: \ifx\theasciiauthors\relax
\theauthors\else\theasciiauthors\fi\s<\ifx\theasciiemail\relax\theemail\else\theasciiemail\fi>}
\immediate\write\gtoutfile{\noexpand\\}
\immediate\write\gtoutfile{Authors: \ifx\theasciiauthors\relax
\theauthors\else\theasciiauthors\fi}
{\def\\{ }\immediate\write\gtoutfile{Title: \ifx\theasciititle\relax
\thetitle\else\theasciititle\fi}}
\immediate\write\gtoutfile{Subj-class: GT or SG, GR etc}
\immediate\write\gtoutfile{MSC-class: \theprimaryclass\ifx\thesecondaryclass\relax\else, \thesecondaryclass\fi}
\immediate\write\gtoutfile{Journal-ref: Algebr. Geom. Topol. \thevolumenumber\s
(\thevolumeyear) \startpage-\finishpage}
\immediate\write\gtoutfile{Comments: Published by Algebraic and
Geometric Topology at}
\immediate\write\gtoutfile{\s\s\s  http://www.maths.warwick.ac.uk/agt/AGTVol\thevolumenumber/agt-\thevolumenumber-\thepapernumber.abs.html}
\immediate\write\gtoutfile{\noexpand\\}
\immediate\write\gtoutfile{}
\ifx\theasciiabstract\relax
\immediate\write\gtoutfile{\theabstract}\else
\immediate\write\gtoutfile{\theasciiabstract}\fi
\immediate\write\gtoutfile{}
\immediate\write\gtoutfile{\noexpand\\}
\immediate\write\gtoutfile{}
\immediate\closeout\gtoutfile}}  

\def\maketitlepage{\makeagttitle\makeheadfile}
\let\makeshorttitle\maketitlepage
\let\maketitle\maketitlepage

\def\ifplaintex{\expandafter\ifx\csname documentclass\endcsname\relax}

\def\gtp{{\mathsurround=0pt\it $\cal G\mskip-2mu$eometry \&\ 
$\cal T\!\!$opology $\cal P\!$ublications}}  

\def\recd{{\small Received:\qua\receiveddate\ifx\reviseddate\relax
\else\qquad Revised:\qua\reviseddate\fi\par}} 

\def\lognumber#1{\def\thelognumber{#1}}
\def\volumenumber#1{\def\thevolumenumber{#1}}
\def\volumeyear#1{\def\thevolumeyear{#1}}
\def\papernumber#1{\def\thepapernumber{#1}}
\def\pagenumbers#1#2{\def\startpage{#1}\def\finishpage{#2}}
\def\published#1{\def\publishdate{#1}}

\def\received#1{\def\receiveddate{#1}}
\def\revised#1{\def\reviseddate{#1}}
\def\accepted#1{\def\accepteddate{#1}}

\def\asciiaddress#1{\def\theasciiaddress{#1}}

\let\\\par\let\thelognumber\relax\let\thevolumenumber\relax
\let\thepapernumber\relax\let\thevolumeyear\relax\let\startpage\relax
\let\finishpage\relax\let\publishdate\relax\let\receiveddate\relax
\let\reviseddate\relax\let\accepteddate\relax\let\theasciititle\relax
\let\theasciiauthors\relax\let\theasciiaddress\relax
\let\theasciiabstract\relax

\let\theasciiemail\relax

\ifplaintex
\font\logobig=cmssbx10 scaled 3836
\font\logomed=cmssbx10 scaled 2557
\else
\font\logobig=cmssbx10 scaled 4200
\font\logomed=cmssbx10 scaled 2800
\fi

\long\def\makeagttitle{   
\count0=\startpage
\agt\hfill      
\hbox to 45truept{\vbox to 0pt{\vglue -13truept{\logomed A\kern -.37em{\logobig 
T}\kern -.38em G}\vss}\hss}
\break
{\small Volume \thevolumenumber\ (\thevolumeyear)
\startpage--\finishpage\nl
Published: \publishdate}

\vglue .25truein

{\parskip=0pt\leftskip 0pt plus
1fil\def\\{\par\smallskip}{\Large\bf\thetitle}\par\medskip} \vglue
0.05truein

{\parskip=0pt\leftskip 0pt plus 1fil\def\\{\par}{\sc\theauthors}
\par\medskip}%
 
\vglue 0.03truein 

{\small\leftskip 25truept\rightskip 25truept{\bf Abstract}\stdspace\theabstract

{\bf AMS Classification}\stdspace\theprimaryclass
\ifx\thesecondaryclass\relax\else; \thesecondaryclass\fi\par
{\bf Keywords}\stdspace \thekeywords\par}\vglue 7truept

}   

\ifplaintex
\font\phead=cmsl9 scaled 950
\font\pnum=cmbx10 scaled 913
\font\pfoot=cmsl9 scaled 950
\headline{\vbox to 0pt{\vskip -4.5mm\line{\small\phead\ifnum
\count0=\startpage ISSN 1472-2739 (on-line) 1472-2747 (printed)
\hfill {\pnum\folio}\else\ifodd\count0\def\\{ }%
\ifx\theshorttitle\relax\thetitle\else\theshorttitle\fi\hfill{\pnum\folio}
\else\def\\{ and }{\pnum\folio}\hfill\ifx\theshortauthors\relax\theauthors
\else\theshortauthors\fi\fi\fi}\vss}}
\footline{\vbox to 0pt{\vglue 0mm\line{\small\pfoot\ifnum\count0=\startpage
\copyright\ \gtp\hfill\else
\agt, Volume \thevolumenumber\ (\thevolumeyear)\hfill\fi}\vss}}
\else
\font=cmsl9 scaled 1050
\font\lnum=cmbx10 
\font=cmsl9 scaled 1050
\makeatletter
\def\@oddhead{{\small\lhead\ifnum\count0=\startpage ISSN 1472-2739 
(on-line) 1472-2747 (printed)\hfill {\lnum\number\count0}\else\ifodd\count0
\def\\{ }\ifx\theshorttitle\relax \thetitle \else\theshorttitle\fi\hfill
{\lnum\number\count0}\else\def\\{ and }{\lnum\number\count0}
\hfill\ifx\theshortauthors\relax 
\theauthors\else\theshortauthors\fi\fi\fi}}\def\@evenhead{\@oddhead}
\def\@oddfoot{\small\lfoot\ifnum\count0=\startpage\copyright\ \gtp\hfill\else
\agt, Volume \thevolumenumber\ (\thevolumeyear)\hfill\fi}
\def\@evenfoot{\@oddfoot}
\makeatother
\fi
\let\maketitlepage\makeagttitle
\let\makeshorttitle\maketitlepage
\let\maketitle\maketitlepage

\newwrite\gtoutfile
\long\gdef\makeheadfile{  
{\def\\{, }\def\s{ }
\immediate\openout\gtoutfile head.xxx
\immediate\write\gtoutfile{To: math@arxiv.org}
\immediate\write\gtoutfile{Subject: put OR rep NNNNN:ppppp}
\immediate\write\gtoutfile{--text follows this line--}
\immediate\write\gtoutfile{Proxy-for: \ifx\theasciiauthors\relax
\theauthors\else\theasciiauthors\fi\s<\ifx\theasciiemail\relax\theemail\else\theasciiemail\fi>}
\immediate\write\gtoutfile{\noexpand\\}
\immediate\write\gtoutfile{Authors: \ifx\theasciiauthors\relax
\theauthors\else\theasciiauthors\fi}
{\def\\{ }\immediate\write\gtoutfile{Title: \ifx\theasciititle\relax
\thetitle\else\theasciititle\fi}}
\immediate\write\gtoutfile{Subj-class: GT or SG, GR etc}
\immediate\write\gtoutfile{MSC-class: \theprimaryclass\ifx\thesecondaryclass\relax\else, \thesecondaryclass\fi}
\immediate\write\gtoutfile{Journal-ref: Algebr. Geom. Topol. \thevolumenumber\s
(\thevolumeyear) \startpage-\finishpage}
\immediate\write\gtoutfile{Comments: Published by Algebraic and
Geometric Topology at}
\immediate\write\gtoutfile{\s\s\s  http://www.maths.warwick.ac.uk/agt/AGTVol\thevolumenumber/agt-\thevolumenumber-\thepapernumber.abs.html}
\immediate\write\gtoutfile{\noexpand\\}
\immediate\write\gtoutfile{}
\ifx\theasciiabstract\relax
\immediate\write\gtoutfile{\theabstract}\else
\immediate\write\gtoutfile{\theasciiabstract}\fi
\immediate\write\gtoutfile{}
\immediate\write\gtoutfile{\noexpand\\}
\immediate\write\gtoutfile{}
\immediate\closeout\gtoutfile}}  

\def\maketitlepage{\makeagttitle\makeheadfile}
\let\makeshorttitle\maketitlepage
\let\maketitle\maketitlepage

\def\ifplaintex{\expandafter\ifx\csname documentclass\endcsname\relax}

\def\gtp{{\mathsurround=0pt\it $\cal G\mskip-2mu$eometry \&\ 
$\cal T\!\!$opology $\cal P\!$ublications}}  

\def\recd{{\small Received:\qua\receiveddate\ifx\reviseddate\relax
\else\qquad Revised:\qua\reviseddate\fi\par}} 

\def\lognumber#1{\def\thelognumber{#1}}
\def\volumenumber#1{\def\thevolumenumber{#1}}
\def\volumeyear#1{\def\thevolumeyear{#1}}
\def\papernumber#1{\def\thepapernumber{#1}}
\def\pagenumbers#1#2{\def\startpage{#1}\def\finishpage{#2}}
\def\published#1{\def\publishdate{#1}}

\def\received#1{\def\receiveddate{#1}}
\def\revised#1{\def\reviseddate{#1}}
\def\accepted#1{\def\accepteddate{#1}}

\def\asciiaddress#1{\def\theasciiaddress{#1}}

\let\\\par\let\thelognumber\relax\let\thevolumenumber\relax
\let\thepapernumber\relax\let\thevolumeyear\relax\let\startpage\relax
\let\finishpage\relax\let\publishdate\relax\let\receiveddate\relax
\let\reviseddate\relax\let\accepteddate\relax\let\theasciititle\relax
\let\theasciiauthors\relax\let\theasciiaddress\relax
\let\theasciiabstract\relax

\let\theasciiemail\relax

\ifplaintex
\font\logobig=cmssbx10 scaled 3836
\font\logomed=cmssbx10 scaled 2557
\else
\font\logobig=cmssbx10 scaled 4200
\font\logomed=cmssbx10 scaled 2800
\fi

\long\def\makeagttitle{   
\count0=\startpage
\agt\hfill      
\hbox to 45truept{\vbox to 0pt{\vglue -13truept{\logomed A\kern -.37em{\logobig 
T}\kern -.38em G}\vss}\hss}
\break
{\small Volume \thevolumenumber\ (\thevolumeyear)
\startpage--\finishpage\nl
Published: \publishdate}

\vglue .25truein

{\parskip=0pt\leftskip 0pt plus
1fil\def\\{\par\smallskip}{\Large\bf\thetitle}\par\medskip} \vglue
0.05truein

{\parskip=0pt\leftskip 0pt plus 1fil\def\\{\par}{\sc\theauthors}
\par\medskip}%
 
\vglue 0.03truein 

{\small\leftskip 25truept\rightskip 25truept{\bf Abstract}\stdspace\theabstract

{\bf AMS Classification}\stdspace\theprimaryclass
\ifx\thesecondaryclass\relax\else; \thesecondaryclass\fi\par
{\bf Keywords}\stdspace \thekeywords\par}\vglue 7truept

}   

\ifplaintex
\font\phead=cmsl9 scaled 950
\font\pnum=cmbx10 scaled 913
\font\pfoot=cmsl9 scaled 950
\headline{\vbox to 0pt{\vskip -4.5mm\line{\small\phead\ifnum
\count0=\startpage ISSN 1472-2739 (on-line) 1472-2747 (printed)
\hfill {\pnum\folio}\else\ifodd\count0\def\\{ }%
\ifx\theshorttitle\relax\thetitle\else\theshorttitle\fi\hfill{\pnum\folio}
\else\def\\{ and }{\pnum\folio}\hfill\ifx\theshortauthors\relax\theauthors
\else\theshortauthors\fi\fi\fi}\vss}}
\footline{\vbox to 0pt{\vglue 0mm\line{\small\pfoot\ifnum\count0=\startpage
\copyright\ \gtp\hfill\else
\agt, Volume \thevolumenumber\ (\thevolumeyear)\hfill\fi}\vss}}
\else
\font=cmsl9 scaled 1050
\font\lnum=cmbx10 
\font=cmsl9 scaled 1050
\makeatletter
\def\@oddhead{{\small\lhead\ifnum\count0=\startpage ISSN 1472-2739 
(on-line) 1472-2747 (printed)\hfill {\lnum\number\count0}\else\ifodd\count0
\def\\{ }\ifx\theshorttitle\relax \thetitle \else\theshorttitle\fi\hfill
{\lnum\number\count0}\else\def\\{ and }{\lnum\number\count0}
\hfill\ifx\theshortauthors\relax 
\theauthors\else\theshortauthors\fi\fi\fi}}\def\@evenhead{\@oddhead}
\def\@oddfoot{\small\lfoot\ifnum\count0=\startpage\copyright\ \gtp\hfill\else
\agt, Volume \thevolumenumber\ (\thevolumeyear)\hfill\fi}
\def\@evenfoot{\@oddfoot}
\makeatother
\fi
\let\maketitlepage\makeagttitle
\let\makeshorttitle\maketitlepage
\let\maketitle\maketitlepage

\newwrite\gtoutfile
\long\gdef\makeheadfile{  
{\def\\{, }\def\s{ }
\immediate\openout\gtoutfile head.xxx
\immediate\write\gtoutfile{To: math@arxiv.org}
\immediate\write\gtoutfile{Subject: put OR rep NNNNN:ppppp}
\immediate\write\gtoutfile{--text follows this line--}
\immediate\write\gtoutfile{Proxy-for: \ifx\theasciiauthors\relax
\theauthors\else\theasciiauthors\fi\s<\ifx\theasciiemail\relax\theemail\else\theasciiemail\fi>}
\immediate\write\gtoutfile{\noexpand\\}
\immediate\write\gtoutfile{Authors: \ifx\theasciiauthors\relax
\theauthors\else\theasciiauthors\fi}
{\def\\{ }\immediate\write\gtoutfile{Title: \ifx\theasciititle\relax
\thetitle\else\theasciititle\fi}}
\immediate\write\gtoutfile{Subj-class: GT or SG, GR etc}
\immediate\write\gtoutfile{MSC-class: \theprimaryclass\ifx\thesecondaryclass\relax\else, \thesecondaryclass\fi}
\immediate\write\gtoutfile{Journal-ref: Algebr. Geom. Topol. \thevolumenumber\s
(\thevolumeyear) \startpage-\finishpage}
\immediate\write\gtoutfile{Comments: Published by Algebraic and
Geometric Topology at}
\immediate\write\gtoutfile{\s\s\s  http://www.maths.warwick.ac.uk/agt/AGTVol\thevolumenumber/agt-\thevolumenumber-\thepapernumber.abs.html}
\immediate\write\gtoutfile{\noexpand\\}
\immediate\write\gtoutfile{}
\ifx\theasciiabstract\relax
\immediate\write\gtoutfile{\theabstract}\else
\immediate\write\gtoutfile{\theasciiabstract}\fi
\immediate\write\gtoutfile{}
\immediate\write\gtoutfile{\noexpand\\}
\immediate\write\gtoutfile{}
\immediate\closeout\gtoutfile}}  

\def\maketitlepage{\makeagttitle\makeheadfile}
\let\makeshorttitle\maketitlepage
\let\maketitle\maketitlepage

\lognumber{18}
\volumenumber{2}
\volumeyear{2002}
\papernumber{18}
\published{22 May 2002}
\pagenumbers{381}{389}
\received{14 January 2002}
\revised{16 May 2002}
\accepted{17 May 2002}

\usepackage{amssymb} 
\usepackage{amscd}
\usepackage[leqno]{amsmath} 
\numberwithin{equation}{section}

\newtheorem{Lemma}[subsection]{Lemma}

\newtheorem{Corollary}[subsection]{Corollary}
\newtheorem{Theorem}[subsection]{Theorem}

\theoremstyle{definition}
\newtheorem{Definition}[subsection]{Definition}
\newtheorem{Remark}[subsection]{Remark}
\newtheorem{Example}[subsection]{Example}

\newcommand{\NE}[1]{\begin{picture}(1.9,1)
\put(2.4,-1){\vector(1,1){2.2}} \put(2.7, 0.15){$\scriptstyle{#1}$}
\end{picture} &&}

\newcommand{\NEE}[1]{\begin{picture}(1.9,1)
\put(.5,-1){\vector(2,1){5}} \put(4, -.5){$\scriptstyle{#1}$}
\end{picture} && }

\newcommand{\So}[1]{\begin{picture}(1.9,1)(0,-0.2) \put(0 ,
1.2){\vector(0,-1){6}} \put(-1, -2){$\scriptstyle{#1}$}
\end{picture}\hspace{-1.9em}}

\newcommand{\SE}[1]{\begin{picture}(1.9,1)(0,0.0) \put(2.4 ,
1.2){\vector(1,-1){2.2}} \put(1.9, -0.05){$\scriptstyle{#1}$}
\end{picture} &&}

\newcommand{\SEE}[1]{\begin{picture}(1.9,1)(0,0.0) \put(.5 ,
1.7){\vector(2,-1){5}} \put(4, .5){$\scriptstyle{#1}$} \end{picture}
&&}

\newcommand\bs{\backslash}

 \newcommand\BR{\mathbb{R}}
 \newcommand\BS{\mathbb{S}}
 
\newcommand\BD{\mathbb{D}}

\def\cor{\begin{Corollary}} \def\defn{\begin{Definition}}
\def\enddefn{\end{Definition}} \def\epr{\end{Theorem}}
\def\ex{\begin{Example}} \def\demo{\begin{proof}}
\def\enddemo{\end{proof}} \def\rem{\begin{Remark}} \def\Ref{\ref}

\DeclareMathOperator{\im}{Im}

\DeclareMathOperator{\GL}{GL} 
\DeclareMathOperator{\s}{S} \renewcommand{\O}{\operatorname{O}}
\DeclareMathOperator{\Top}{Top} \DeclareMathOperator{\Int}{Int}

\begin{document}

\title{Exotic smooth structures on nonpositively\\curved symmetric
spaces} 
\author{Boris Okun} 

\address{Department of Mathematical
Sciences, University of Wisconsin--Milwaukee\\Milwaukee, WI 53201, USA}

\asciiaddress{Department of Mathematical
Sciences, University of Wisconsin-Milwaukee\\Milwaukee, WI 53201, USA}

\email{okun@uwm.edu}

\begin{abstract} We construct series of examples of exotic smooth
structures on compact locally symmetric spaces of noncompact type. In
particular, we obtain higher rank examples, which do not support
Riemannian metric of nonpositive curvature. The examples are obtained
by taking the connected sum with an exotic sphere. To detect the change
of the smooth structure we use a tangential map from the locally
symmetric space its dual compact type twin.  \end{abstract}



\maketitle

\section{ Introduction}

The purpose of this note is to show how tangential maps between dual
symmetric spaces, constructed in \cite{O1}, can be used to obtain
exotic smooth structures on a compact locally symmetric space of
noncompact type. In particular, we construct higher rank examples,
which then by Eberlein--Gromov Rigidity Theorem do not support a
Riemannian metric of nonpositive curvature. This answers in the
negative the question, due to Eberlein, of whether a smooth closed
manifold, homotopy equivalent to a nonpositively curved one, admits a
Riemannian metric of nonpositive curvature. Rank $1$ examples with the
same property were independently constructed in \cite{AF}.

The simplest way to attempt to change the smooth structure on a
manifold is to take the connected sum with an exotic sphere. The
question remains, however, of whether this procedure actually {\em
changes} the structure. This is where tangential maps become useful.

The paper is organized as follows. Section \ref{TMSS} contains general
statements, showing how to produce nontrivial (up to concordance)
smooth structures of the above sort on the domain of a tangential map
from the ones on its range. In Section \ref {APP}, we apply these
techniques to nonpositively curved locally symmetric spaces and obtain
nondiffeomorphic structures.

Though our arguments are straightforward generalizations of some of the
arguments in \cite{FJ1}, \cite{FJ3}, and are probably well known to
specialists, for reader's convenience we present them in some detail.

The material in this paper is part of the author's 1994 doctoral
dissertation written at SUNY Binghamton under the direction of F.\,T.
Farrell. I wish to thank him and S.\,C. Ferry for numerous useful
discussions. The author was partially supported by NSF grant.

\section{\label{TMSS} Tangential maps and smooth structures}

In this section, we show that certain types of smooth structures remain
nontrivial under tangential maps. The arguments here are very close to
those of \cite{FJ3}, though our statements are more general.

First, we recall the following very useful lemma:

\begin{Lemma}[{\cite[p. 509]{KM}}] \label{Milnor} Let $\xi$ be a
$k$-dimensional vector bundle over an $n$-dimensional space, $k>n$.  If
the Whitney sum of $\xi$ with a trivial bundle is trivial, then $\xi$
itself is trivial.  \end{Lemma}

Let $\BD^n$ denote the closed ball of radius 1 in the Euclidean
$n$-dimensional space.

\begin{Lemma}{\label{emb}} Let $k\co M^n \to N^n$ be a tangential map
between two closed smooth $n$-dimensional manifolds. Then $M^n \times
\BD^{n+1}$ is diffeomorphic to a codimension $0$ submanifold of the
interior of $N^n \times \BD^{n+1}$.  \end{Lemma}

\demo Let $i\co N^n \to N^n \times \BD^{n+1}$ be the standard inclusion
$i(x)=(x,0)$ and $p\co N^n \times \BD^{n+1} \to N^n$ be the projection
on the first factor $p(x,y)=x$.

Consider the composition $i\circ k\co M^n \to N^n \times \BD^{n+1}$.
By Whitney Embedding Theorem this composition can be approximated by an
embedding $w\co M^n \to N^n \times \BD^{n+1}$, such that $w$  is
homotopic to $i\circ k$.

Let $\nu$ denote the normal bundle of the manifold $M^n$ considered as
a submanifold of $ N^n \times \BD^{n+1}$ via $w$. By definition of the
normal bundle, we have: $$ w^*(\tau( N^n \times \BD^{n+1}))=\nu \oplus
\tau(M^n) $$

On the other hand, $$ \tau( N^n \times \BD^{n+1})=p^*(\tau(N^n)) \oplus
\epsilon^{n+1} $$ where $\epsilon^{n+1}$ denotes the
$(n+1)$-dimensional trivial bundle.

Combining these two equations together, we obtain:  $$ \nu \oplus
\tau(M^n)=  w^*(p^*(\tau(N^n))) \oplus \epsilon^{n+1}$$ Since by
construction $w \simeq i\circ k$, we have $w^*$=$(i\circ k)^*=k^* \circ
i^*$. Note that $ i^* \circ p^* = (p \circ i)^*=id$ and $k^*
\tau(N^n)=\tau(M^n)$ since the map $k$ is tangential. It follows that
$$ \nu \oplus \tau(M^n)= \tau(M^n) \oplus  \epsilon^{n+1}$$ i.e. the
bundle $\nu$ is stably trivial.

By Lemma \Ref{Milnor} $\nu$ is trivial itself; therefore, the tubular
neighborhood of $M^n$ in $N^n \times \BD^{n+1}$ is diffeomorphic to
$M^n \times \BD^{n+1}$.  \enddemo

Let $Z^n$ be a closed, connected $n$-dimensional manifold. There is a
unique homotopy class of degree $1$ maps from $Z^n$ to $\BS^n$, which
we will denote by $f_Z$. One can obtain such a map by collapsing to a
point the exterior of the small neighborhood of a point in $Z^n$.

It turns out that a tangential map induces a map on suspensions going
in the ``wrong'' direction.

\begin{Lemma}[\label{Backwards}cf. {\cite[Corollary 3.11]{FJ3}}] Let
$k\co M^n \to N^n$ be a tangential map between $n$-dimensional
manifolds. Then there exist a map $$g\co \Sigma^{n+1}N^n \to
\Sigma^{n+1}M^n $$ such that the suspension $\Sigma^{n+1}f_N$ and the
composite $(\Sigma^{n+1}f_M)\circ g$ are homotopic as maps from
$\Sigma^{n+1}N^n$ to $\Sigma^{n+1}\BS^n=\BS^{2n+1}$.  \end{Lemma}

\demo Fix base points $x_0 \in M^n$ and $y_0 \in N^n$. Let $B$ denote a
small neighborhood of $x_0$ in $M^n$. By general position, we may
assume that the embedding $$F\co M^n \times \BD^{n+1} \to N^n \times
\BD^{n+1}$$ of Lemma \Ref{emb} has the additional property that:
$$\im(F)\cap (y_0 \times \BD^{n+1}) \subseteq F(B \times \BD^{n+1})$$

The suspensions $\Sigma^{n+1}M^n$ and $\Sigma^{n+1}N^n$ can be
identified as the following quotient spaces:  $$\Sigma^{n+1}M^n=M^n
\times \BD^{n+1}/(M^n \times \partial \BD^{n+1}\cup B \times
\BD^{n+1})$$

$$\Sigma^{n+1}N^n=M^n \times \BD^{n+1}/(M^n \times \partial
\BD^{n+1}\cup y_0 \times \BD^{n+1}) \label{sus} $$

Let $*$ denote the point in $\Sigma^{n+1}M^n$ corresponding to the
subset $M^n \times \partial \BD^{n+1}\cup B \times \BD^{n+1}$ in the
first formula. Define  $g\co \Sigma^{n+1}N^n \to \Sigma^{n+1}M^n $ by
$$g(y)= \begin{cases} F^{-1}(y), &\text{if } y\in F( (M^n-B) \times
\Int (\BD^{n+1} ))\\
	       *, &\text{otherwise.} \\
\end{cases} $$

It is easy to check that $\Sigma^{n+1}f_N$ and $(\Sigma^{n+1}f_M)\circ
g$ are homotopic.  \enddemo

Now let us turn our attention to smooth structures.

\begin{Definition} Let $Z$ be a compact closed topological manifold,
and let $Z_1$, $Z_2$ be two smooth structures on $Z$. $Z_1$ is said to
be {\em concordant} to $Z_2$ if there is a smooth structure $\bar Z$ on
$Z\times [0,1]$ such that $\partial_- \bar Z=Z_1$ and $\partial_+ \bar
Z=Z_2$.  \enddefn

Let $\Sigma^n$ be an exotic $n$-dimensional sphere, \cite{KM}. One of
the ways to produce a (possibly) nontrivial smooth structure on a
smooth $n$-dimensional manifold $Z^n$ is to take the connected sum
$Z^n\#\Sigma^n$. The following lemma shows that, at least for this type
of variation of smooth structures, tangential maps preserve
nontriviality of smooth structures.

\begin{Lemma}{\label{Dual}} Let $k\co M^n \to N^n$ be a tangential map
between n-dimensional manifolds and assume $n \geq 7$. Let $\Sigma_1$
and $\Sigma_2$ be homotopy $n$-spheres. Suppose that $M^n\#\Sigma_1^n$
is concordant to $M^n\#\Sigma_2^n$, then $N^n\#\Sigma_1^n$ is
concordant to $N^n\#\Sigma_2^n$.  \end{Lemma}

We note that there is {\em no condition} on the degree of the map $k$
above.

\demo Recall \cite[p. 194]{KS} that concordance classes of smooth
structures on a smooth manifold $Z^n$ ($n>4$) are in one-to-one
correspondence with homotopy classes of maps from $Z^n$ to $\Top/\O$
denoted by $[Z^n, \Top/\O]$. The space $\Top/\O$ is an infinite loop
space (see \cite[p. 215]{BV}). In particular, $\Top/\O=\Omega^{n+1}Y$
for some space $Y$. We have suspension isomorphism $[Z^n, \Top/\O]
\cong  [\Sigma^{n+1}Z^n,Y]$, so the smooth structures on $Z^n$ can be
thought of as (homotopy classes of) maps from $\Sigma^{n+1}Z^n$ to
$Y$.  In this way a connected sum $Z^n\#\Sigma^n$ gives rise to a map
from $\Sigma^{n+1}Z^n$ to $Y$, which we will denote $z_\Sigma$. Using
the standard sphere $\BS^n$ in place of $Z^n$ we see that an exotic
sphere $\Sigma^n$ itself corresponds to a map $s_\Sigma\co \BS^{2n+1}
\to Y$.  The naturality of this construction and Lemma \Ref{Backwards}
imply that the diagram $$ \begin{CD} \Sigma^{n+1}N^n \\ \So{g}
\SE{n_{\Sigma_i}} \SEE{\Sigma^{n+1}f_N} \\ && Y\hspace{-2em}
@<s_{\Sigma_i}<< \BS^{2n+1} \\ \NE{\hspace{-1.2em} m_{\Sigma_i}}
\NEE{\Sigma^{n+1}f_M} \\ \Sigma^{n+1}M^n \end{CD} $$

is homotopy
commutative for $i=1,2$. In particular, we see that $n_{\Sigma_i}
\simeq m_{\Sigma_i} \circ g$.

Since the connected sums $M^n\#\Sigma_1^n$ and $M^n\#\Sigma_2^n$ are
concordant, the maps $m_{\Sigma_1}$ and $m_{\Sigma_2}$ are homotopic.
It follows that the maps $n_{\Sigma_1}$ and $n_{\Sigma_2}$ are
homotopic, and therefore the connected sums $N^n\#\Sigma_1^n$ and
$N^n\#\Sigma_2^n$ are concordant.  \enddemo

\section{\label{APP} Exotic smoothings of locally symmetric spaces}

In this section, we apply the results of the previous section to the
tangential map constructed in \cite{O1}. We show how to produce
nonstandard smooth structures on the nonpositively curved symmetric
spaces.

We will use notation from \cite{O1}. Let $G$ be a real semisimple
algebraic linear Lie group (a subgroup of $\GL(n,\BR)$) and  $K$ be its
maximal compact subgroup. Let $G_c$ denote the complexification of $G$
and let $G_u$ denote a maximal compact subgroup of $G_c$.  Let $\Gamma$
be a cocompact torsion-free discrete subgroup of $G$. The space
$X=\Gamma \bs G/K$ is a locally symmetric space of noncompact type and
the space $X_u=G_u/K$ is a symmetric space of compact type.  Slightly
abusing terminology, we will refer to the spaces $X$ and $X_u$ as dual
symmetric spaces.

The tangential map is provided by the following theorem:

\begin{Theorem}{\rm\cite{O1}\label{ta}}\qua Let $X=\Gamma\bs G/K$ and
$X_u=G_u/K$ be dual symmetric spaces. Then there exist a finite sheeted
cover $X'$ of $X$ (i.e. there is a subgroup $\Gamma'$ of finite index
in $\Gamma$ with $X'=\Gamma' \bs G/K$) and a tangential map $k\co X'
\to X_u$.  \end{Theorem}

We will also use the following three rigidity results.

\begin{Theorem}[The Strong Mostow Rigidity Theorem \cite{Mo}] Let $X_1$
and $X_2$ be compact locally symmetric spaces of noncompact type such
that the universal cover of $X_1$ has no $2$-dimensional metric factor
projecting to a closed subset of $X_1$. Then any isomorphism from
$\pi_1(X_1)$ to $\pi_1(X_2)$ is induced by a unique isometry (after
adjusting the normalizing constants for $X_1$).  \end{Theorem}

\begin{Theorem}[Eberlein--Gromov Strong Rigidity Theorem 
\cite{E1, BGS}\label{EG}] Let $X$ be compact locally symmetric space
of noncompact type such that all metric factors of $X$ have rank
greater than $1$. Let $M$ be a closed connected nonpositively curved
Riemannian manifold. Then any isomorphism from $\pi_1(X)$ to
$\pi_1(M)$ is induced by a unique isometry (after adjusting the
normalizing constants for $X$).  \end{Theorem}

\begin{Theorem}[Farrell--Jones Topological Rigidity Theorem \cite{FJ2}]
Let $M$ be a closed connected $m$-dimensional Riemannian manifold with
nonpositive sectional curvature, let $N$ be a topological manifold
(possibly with boundary) and let $H\co N \to M\times \BD^n$  be a
homotopy equivalence  which is a homeomorphism on the boundary. Assume
$m+n \neq 3,4$. Then $H$ is homotopic rel boundary to a homeomorphism
from $N$ to $M\times \BD^n$.  i \end{Theorem}

The following lemma provides a connection between concordance and
diffeomorphism classes of smooth structures for locally symmetric
spaces.

\begin{Lemma}[\label{Rigid}cf. \cite{FJ1}] Let $X=\Gamma \bs G/K$ be a
compact orientable symmetric space of noncompact type such that the
universal cover $G/K$ of $X$ has no $2$-dimensional metric factor
projecting to a closed subset of $X$ and assume $\dim X\geq 7$. Let
$\Sigma_1$ and $\Sigma_2$ be homotopy spheres of the same dimension as
$X$. Suppose $X\#\Sigma_1$ is diffeomorphic to $X\#\Sigma_2$. Then
$X\#\Sigma_1$ is concordant either to $X\#\Sigma_2$ or to
$X\#(-\Sigma_2)$.  \end{Lemma}

\demo The argument given to prove Addendum 2.3 in \cite{FJ1} holds
almost literally in our setting. For completeness we give the argument
here.

We can assume that the connected sums of $X$ with $\Sigma_1$ and
$\Sigma_2$ are taking place on the boundary of a small metric ball $B$
in $X$, so we think of  $X\#\Sigma_1$ and $X\#\Sigma_2$ as being
topologically identified with $X$, and the changes in the smooth
structure happen inside $B$. Let $f\co X\#\Sigma_1 \to X\#\Sigma_2$ be
a diffeomorphism.

First we consider a special case, where $f\co X \to X$ is homotopic to
the identity. Thus we have a homotopy $h \co X\times [0,1] \to X$ with
$h\vert_{X\times 1}=f$ and $h\vert_{X\times 0}=id_X$. Define $H \co
X\times [0,1] \to X\times [0,1]$ by $H(x,t)=(h(x,t),t)$. Note that $H$
is a homotopy equivalence which restricts to a homeomorphism on the
boundary. Therefore, by Farrell--Jones Topological Rigidity Theorem,
$H$ is homotopic rel boundary to a homeomorphism $H' \co X\times [0,1]
\to X\times [0,1]$. We put the smooth structure $(X\#\Sigma_2) \times
[0,1]$ on the range of $H'$ and, by pulling it back along $H'$, obtain
the smooth structure $N$ on the domain of $H'$. Since, by
construction,  $H' \co N \to (X\#\Sigma_2) \times [0,1]$ is a
diffeomorphism, $H'\vert_{X\times 1}=f$ and $H'\vert_{X\times 0}=id_X$,
$N$ is a concordance between $X\#\Sigma_1$ and $X\#\Sigma_2$.

The general case reduces to the above special case as follows. If $f
\co X \to X$ is an orientation preserving homeomorphism, then, by the
Strong Mostow Rigidity Theorem, $f^{-1}$ is homotopic to an orientation
preserving isometry $g$. Since one can move around small metric balls
in $X$ by smooth isotopies, $g$ is homotopic to a diffeomorphism $g'
\co X \to X$ such that $g'\vert_B=id_B$. Since $X\#\Sigma_2$ is
obtained by taking the connected sum along the boundary of $B$ and
$g'\vert_B=id_B$ and $g'\vert_{M-B}$ is a diffeomorphism, it follows
that $g' \co X\#\Sigma_2  \to X\#\Sigma_2$  is also a diffeomorphism.
Therefore the composition $ g'\circ f \co X\#\Sigma_1 \to X\#\Sigma_2$
is a diffeomorphism homotopic to the identity and it follows from the
previous special case that $X\#\Sigma_1$ and $X\#\Sigma_2$ are
concordant.

If $f$ is an orientation reversing homeomorphism, then similar argument
produces a diffeomorphism $g' \co X\#\Sigma_2 \to X\#(-\Sigma_2)$  and
it follows that $X\#\Sigma_1$ is concordant to $X\#(-\Sigma_2)$.
\enddemo

Our main result is the following theorem:

\begin{Theorem}{\label{Smooth}} Let $X=\Gamma\bs G/K$ and $X_u=G_u/K$
be compact dual symmetric spaces such that the universal cover $G/K$ of
$X$ has no $2$-dimensional metric factor projecting to a closed subset
of $X$ and assume $\dim X\geq 7$. Let $X'$ be the oriented finite
sheeted cover of $X$ the existence of which was established by Theorem
\Ref{ta}. Let $\Sigma_1$ and $\Sigma_2$ be homotopy spheres of the same
dimension as~$X$.

If the connected sum $X_u\#\Sigma_1$ is not concordant to both
$X_u\#\Sigma_2$ and $X_u\#(-\Sigma_2)$ then $X'\#\Sigma_1$ and
$X'\#\Sigma_2$ are not diffeomorphic.  \end{Theorem}

\demo Suppose the connected sums $X'\#\Sigma_1$ and $X'\#\Sigma_2$ are
diffeomorphic. Then, by Lemma \Ref{Rigid}, $X'\#\Sigma_1$ is concordant
either to $X'\#\Sigma_2$ or to $X'\#(-\Sigma_2)$. It follows from Lemma
\Ref{Dual} that either $X_u\#\Sigma_1$ is concordant to $X_u\#\Sigma_2$
or to $X_u\#(-\Sigma_2)$, which contradicts the hypothesis. This
contradiction proves the theorem.  \enddemo

Applying this theorem to \cite[Example 4]{O1}, we obtain the following
corollary:

\cor{\label{example}} Let $G=G_c$ be a complex semisimple Lie group,
and let $X=\Gamma\bs G_c/G_u$ and $X_u=G_u$ be compact dual symmetric
spaces such that $\dim X\geq 7$. Let $X'$ be the oriented finite
sheeted cover of $X$ the existence of which was established by Theorem
\Ref{ta}. Let $\Sigma_1$ and $\Sigma_2$ be two nondiffeomorphic
homotopy spheres of the same dimension as $X$. Then the connected sums
$X'\#\Sigma_1$ and $X'\#\Sigma_2$ are not diffeomorphic.
\end{Corollary}

\demo Since the group $G_c$ is a complex Lie group, we do not have to
worry about $2$-dimensional metric factors; hence, according to Theorem
\Ref{Smooth} it suffices to show that the connected sum $G_u\#\Sigma_1$
is not concordant to both $G_u\#\Sigma_2$ and $G_u\#(-\Sigma_2)$. We
recall that the tangent bundle of a Lie group is trivial. Therefore the
constant map $G_u \to \BS^{\dim G_u}$ is tangential, and the required
statement follows from Lemma \Ref{Dual}.  \enddemo

The above corollary provides examples of exotic smooth structures on
locally symmetric space of higher rank, therefore by Eberlein--Gromov
Rigidity Theorem \ref{EG} these examples do not admit Riemannian metric
of nonpositive curvature. One can apply Theorem \Ref{Smooth} to real
and complex hyperbolic manifolds to obtain exotic smoothings of
Farrell--Jones \cite{FJ1},\cite{FJ3}. This is trivial in the real case,
since the dual space is a sphere, and requires rather intricate
analysis of the smooth structures on the complex projective space (see
\cite{FJ3}) in the complex case. The case of the quaternionic
hyperbolic manifolds seems still to be open.


\Addresses\recd
\end{document}